\documentclass[11pt, reqno]{amsart}
\usepackage{amscd}
\usepackage[enableskew,vcentermath]{youngtab}
\usepackage{amssymb}
\usepackage{amsmath}
\usepackage[all]{xypic}
\usepackage{graphics}

\def\eee{\mbox{$\Box$}}

\def\ev{{\mathchoice{\mbox{\rm ev}}
                    {\mbox{\rm ev}}
                    {\mbox{\scriptsize\rm ev}}
                    {\mbox{\tiny\rm ev}} }}

\def\id{{\mathchoice{\mbox{\rm id}}
                    {\mbox{\rm id}}
                    {\mbox{\scriptsize\rm id}}
                    {\mbox{\tiny\rm id}} }}

\newtheorem{thm}{Theorem}[section]

\newtheorem{dn}[thm]{Definition}

\newtheorem{pro}[thm]{Proposition}
\newcommand{\bbar}{\begin{array}}
\newcommand{\eear}{\end{array}}
\newcommand{\bb}{\begin{equation}}
\newcommand{\eqbb}{\begin{equation}}
\def\ee{\end{equation}}
\def\eqee{\end{equation}}
\def\eea{\end{eqnarray}}
\def\bba{\begin{eqnarray}}
\def\Ker{\mbox{\rm Ker}}

\def\det{\mbox{\rm det}}
\def\ch{\mbox{\rm ch}}

\def\Im{\mbox{\rm Im}}

\def\End{\mbox{\rm End\hspace{0.1ex}}}

\def\Id{\id}

\def\secteqn{
\let\sectio\section%
\renewcommand{\section}{\sectioneqn\sectio }%
\newcommand{\sectioneqn}{\setcounter{equation}{0}%
 \renewcommand{\theequation}{\arabic{section}.\arabic{equation}}}}\oddsidemargin0.5cm
\title[Double Koszul Complex]{Double Koszul Complex and Construction of Irreducible Representations of   $\frak{gl}(3|1)$ }

\author[N.T.P.Dung]{NGUY$\tilde{\mbox{\^E}}$N Thi
Phuong Dung}

 \address[N.T.P.Dung]{Border Academy, Son Tay, Ha Noi, Vietnam. 
Email: phuongdung72@yahoo.com}

\thanks{Financial support provided to the author by NAFOSTED under gand number 101.01.16.09}

\begin{document}
 \maketitle
 \bibliographystyle{plain}
\section{Introduction}
Let $V$ be a super vector space over a field $k$ of characteristic of $0$. The super group $GL(V)$ of linear  automorphisms of $V$ is the 
subgroup of the semi-group $\End(V)$ of endomorphisms with invertible super-determinant. In \cite{Manin1}, Manin introduced  the following Koszul 
complex $K$ to define the super determinant. 
Its $(k,l)$-term is given by \break $K^{k,l} := \Lambda _k\otimes S_l^*$, where $\Lambda_n$ and  $S_n$ are the $n$-th homogeneous components of
the exterior and the symmetric tensor algebras on $V$. The differential $d_{k,l}: K^{k,l} \longrightarrow K^{k+1,l+1}$ is given by 
$$
d_{k,l} (h\otimes \varphi ) = \sum_i h\wedge x_i \otimes \xi ^i\cdot \varphi .
$$

There is another Koszul complex associated to $V$, denoted by $L$. 
This complex was first defined by Priddy as a free resolution of $k$ as a module over the symmetric tensor
algebra of $V$, see \cite{Manin}. 
  Its $(l,k)$-term is given by  $L^{l,k} := S_l\otimes \Lambda_k$ with  differential $P_{l,k}: L^{l,k} \longrightarrow L^{l-1,k+1}$ given by
$$ \xymatrix { P_{l,k}: S_l\otimes  \Lambda_k \ar@{^(->}[r]&V^{\otimes l}\otimes  V^{\otimes k}=V^{\otimes l-1}\otimes  V^{\otimes k+1}
\ar[rr]^{\qquad \quad X_{l-1}\otimes Y_{k+1}}&&
S_{l-1}\otimes \Lambda_{k+1} },$$
 where $ X_l, Y_k$ are the symmetrizer and anti-symmetrizer operators.
  In \cite{Kac3}, Kac  proved that any finite dimensional irreducible representation  of  the Lie super algebra
  $\frak{gl}(V)$  is a quotient of the Kac module. He divided  irreducible representations of $\frak{gl}(V)$ into two classes, 
  typical representations and atypical representations.
   By using the  Kac module, Kac   gave an explicit construction of all typical representations of $\frak{gl}(V)$, and
a character formula for all typical representations.  
 In  \cite{Zhang2},  Su and   Zhang gave a character formula for all finite-dimensional irreducible representations of  $\frak{gl}(V)$.   An  explicit construction of atypical
  representations is however not known.
 The aim of this work is to give a combinatorial way to describe
all irreducible representations  in case the super-dimension of $V$ is $(3|1)$.

 Our observation   is   that the  two Koszul  complexes above  can be combined into a 
double complex which we call  the double Koszul complex.  We use the differential of this complex to describe all irreducible representations of
$\frak {gl}(V)$ when $V$ has super-dimension $(3|1)$.

The paper is organized as follows.  Section $2$ provides some background materials on the general linear super-algebra needed for the rest of the paper.
Section 3 introduces and studies the double Koszul complex.  Section 4 uses the properties of the double Koszul complex to construct representations of the
Lie super algebra  $\frak {gl}(V)$. Using  the character formula of  Su and Zhang in \cite{Zhang2}, we prove that the constructed representations furnish all
irreducible representations of $\frak {gl}(V)$.

{\it Acknowledgements:} I would like to thank  my teacher, Prof  Phung Ho Hai for his  constant help and suggestions.  I also would like to thank the refree for  usefull comments.

\section{Preliminaries}  
 This section presents some results on the general linear Lie super-algebras for  the later use. We shall work with a field $k$ of characteristic $0$.
  A  super vector space is  a \break $\mathbb {Z}_2$-graded vector space $V =  V_{\bar 0}\oplus V_{\bar 1} .$
  The spaces $ V_{\bar 0},   V_{\bar 1}$ 
 are called  even and odd homogeneous components of $V$, their elements are
 called homogeneous. We   denote the $\mathbb Z_2$-grade (or parity) of a homogeneous element $a$ by $\hat a$.
Assume  $\dim V_{\bar 0}=m, \dim V_{\bar 1}=n$ and  fix a homogeneous basis of $V$: $x_1,\ldots,x_m \in 
 V_{\bar 0}$, $ x_{m+1},\ldots,x_{m+n}\in V_{\bar 1}.$ For simplicity we  denote the $\mathbb Z_2$-grade of $x_i$ by
 $\hat i$. Thus $\hat i=\bar 0$ if $1\leq i\leq m$ and $\hat i=\bar 1$ if $m+1\leq i\leq m+n$.

 A $\mathbb Z_2$-graded algebra $A$ is called a super algebra.
 Similarly we have the notion of super Lie algebra $L$, where the super anti-commutativity and the super Leibniz rule read:
 $$\begin{array}{rcl} [a,b]&=&(-1)^{\hat a\hat b}[b,a],\\[0em]
 [a,[b,c]]&=&[[a,b],c]+(-1)^{\hat a\hat b}[b,[a,c]].\end{array}$$
 Here we use the convention that $(-1)^{\bar 0}=1$ and $(-1)^{\bar 1}=-1$.
 
 Given a super algebra $A$, the super-commutator on $A$, defined by
 $$[a,b]:=ab-(-1)^{\hat a \hat b}ba,$$
 makes $A$ into a super Lie algebra, denoted by $A^L$.
\subsection{Super Lie Algebras $\frak g=\frak{gl}(V)$}   Consider the algebra $\mathrm{End}(V)$
of linear endomorphisms of $V$.
 Fix a homogeneous basis of $V$ as above. Every element of $\mathrm{End}(V)$ is given by a matrix
of the form $ \left({A\ B\atop C\ D}\right)$, where $A,B,C,D$ are block matrices. The matrices of the form 
$ \left({A\ 0\atop 0\ D}\right)$ define  even maps $V\to V$ (i.e.  maps that preserve  the $\mathbb Z_2$-grading).
The matrices of the form $ \left({0\ B\atop C\ 0}\right)$ define odd maps (i.e. maps that interchange the $\mathbb Z_2$-grading).
An arbitrary map $V\to V$ is the sum of an even map  with an odd map. This
defines a $\mathbb Z_2$-grading on $\mathrm{End}(V)$  and makes $\mathrm{End}(V)$ a super  algebra.
The associated super Lie algebra $\mathrm{End}(V)^L$ is denoted by $\frak{gl}(V)$.
\subsection{ Representation of  $\frak g= \frak {gl}(V)$}
  Let $W$ be a super vector space. A super representation $\rho $ of  $\frak g$ in $W$
is an even linear mapping $\rho : \frak g\longrightarrow \frak{gl}(W)$  which preserves the super commutator,
that is a homomorphism of Lie super algebras. A super representation of $\frak g$ is also called a $\frak g$-module. A super representation is said to be irreducible if it has no proper
non-zero sub-representations. In oder to construct all irreducible representations of $\frak g$ we  need the
technique of induced representations, which we will now describe.  
\subsubsection{Induced representations}
A pair $(\mathcal{U}(\frak g), i)$, where 
$ \mathcal{U}(\frak g)$ is an associative $\mathbb{Z}_2$-graded algebra and $i: \frak g \longrightarrow \mathcal{U}(\frak g)^L$ is a 
homomorphism of Lie super algebra, is called a universal enveloping super algebra of $\frak g$ if for any other pair $(\mathcal{U}',i')$, there is an unique
 homomorphism  $\theta : \mathcal{U}\longrightarrow \mathcal{U}'$ such that $i' = \theta .i$. Thus, the concepts  of ``super representation of $\frak g$'', 
 "$\frak g$-module'' and  "left  $\mathcal{U}(\frak g)$-module'' are completely equivalent. 

 Let  $\frak g$ be a super Lie algebra, $\mathcal{U}(\frak g)$ be its universal enveloping super algebra, $\frak h$ be a super Lie sub-algebra of $ \frak g,$ 
and $V$ be
 an $\frak h$-module. The $\mathbb{Z}_2$-graded space $\mathcal{U}(\frak g)\otimes _{\mathcal{U}(\frak h)} V$ can be endowed 
with the structure of a $\frak g$-module as follows: $g (u\otimes v) = g u \otimes v$ for $ g \in \frak g, u\in \mathcal{U}(\frak g), v \in V.$
The so constructed $\frak g$-module is said to be induced from the $\frak h$-module $V$ and is denoted by $\mathrm{Ind}_{\frak h}^{\frak g} V$.
\subsubsection{Weights and Roots of $\frak g$.}
The standard basis for $\frak g$ consists of matrices $E_{ij}: i,j = 1\ldots, m+n$ where $E_{ij}$ is the matrix with $1$
in the place $(i,j)$ and $0$ elsewhere.
 Consider the sub-algebra $\frak h$
of $\frak g$   spanned by the elements $h_j := E_{jj}: i,j = 1, \ldots ,m+n$, $\frak h$ is a Cartan subalgebra of $\frak g$.  
The  space $\frak h^*$  dual to  $\mathcal {\frak h}$ is  is spanned by $\epsilon _i : i  = 1, \ldots ,m+n,$ where  for $X=\left({A\ B\atop C\ D}\right)$,
$$\begin{array}{l}\epsilon _i: X\longmapsto A_{ii},\text{ for } 1\leq i\leq m, \quad
\epsilon _j:X\longmapsto D_{jj},  \text{ for }  m+1\leq j\leq   m+n.\end{array}$$

Elements of $\mathcal {\frak h}^*$ are called the weights of $\frak g$.
Let $\lambda \in \frak h^*$,  
$\lambda =\sum_{i+1}^m \lambda _i\epsilon _i - \sum_{j= m+1}^{m+n}\lambda _j\epsilon _j$
then we write
$\lambda = ( \lambda _1,\lambda _2,\ldots, \lambda _m| \lambda _{m+1}, \ldots ,\lambda _{m+n}).$
 
\begin{dn}
Let $\lambda = (\lambda _1,\lambda _2,\ldots, \lambda _m| \lambda _{m+1},\ldots ,\lambda _{m+n})$ be a weight.
\begin{itemize}\item[(i)] $\lambda $ is called integral   if 
$\lambda _i - \lambda _{i+1} \in \mathbb{Z}$ for  all $i\neq m$. 
\item[(ii)] $\lambda $ is called dominant   if $\lambda _i \geq \lambda _{i+1}$ for  $ 1\leq i\leq m$, and 
$\lambda _j\leq \lambda _{j+1}$ for   $m+1\leq j\leq m+n-1$. 
\item[(iii)]
 $\lambda  $  is called typical     if  $ (\lambda _i + m+1-i) - (\lambda _{m+p} + p) \neq 0$  for  all $1\leq i\leq m, 1\leq p\leq n,$ otherwise it is called atypical.
 \item[(iiii)]
 $\lambda$ is called integrable    if $\lambda_i \in \mathbb Z$ for all $i$. 
  \end{itemize}
 \end{dn}

Let  $0\neq  \alpha \in \frak h^*$. Set 
$ \frak g_{\alpha } := \{ a \in \frak g: [h,a] = \alpha (h) a, \forall  h \in \frak h \}.$
  If $\frak g_{\alpha }\neq 0$, then $\alpha$ has 
 the form $\epsilon_i - \epsilon _j : i \neq j$. It is called a root. 
   We set $\Delta _0^+ = \{ \epsilon _i - \epsilon _j : 1\leq i<j\leq m \quad  \mbox{or} \quad m+1 \leq i<j\leq m+n \}, 
\Delta _1^+ = \{ \epsilon _i - \epsilon _j : 1\leq i\leq m \quad  \mbox{and } \quad m+1 \leq  j \leq m+n \} \text { and }$  
$$\textstyle\rho := (m,m-1,\ldots,1|1,2,\ldots,n) - \frac{m+n+1}{2}(1,1,\ldots,1|1,1,\ldots,1).$$
\subsubsection {Kac module}
  For every integral dominant weight $\lambda$, we denote by $V^0(\lambda )$ the finite dimension irreducible $\frak g_{\bar 0}$-module with highest
 weight $\lambda $, $ V^0(\lambda )$  is the $(\frak g_{\bar 0} \oplus \frak g_{+1})$- module  with $\frak g_{+1}$ acting by 0, 
where $  \frak g_{+1}$ is the set of matrices of the form 
$\left({0\ B\atop 0\ 0}\right).$
Set $\bar V(\lambda ) := \mathrm{Ind}^{\frak g}_{\frak g_{\bar 0} + \frak g_{+1}} V^0(\lambda )$. $\bar V(\lambda ) $ 
 contains a unique maximal submodule $M(\lambda ),$ and we set 
$$ V(\lambda ) : = \bar V(\lambda )/ M(\lambda ).$$
 Then $V(\lambda)$ is an  irreducible representation
 with highest weight $\lambda$. The module $  \bar V(\lambda )$   is called  generalized Verma module or Kac module 
 \cite{Kac3}. Kac showed that the $V(\lambda)$'s furnish all
 irreducible $\frak g$-modules of finite dimension.

If $\lambda $  is typical weight then 
$M_\lambda=0$, thus 
$V(\lambda )=\bar V(\lambda)$, in this case $V(\lambda)$ is called typical.
On the other hand, if $\lambda$ is atypical, an explicit construction of $M(\lambda)$ is
not known.  

  \subsubsection{Characters of representations}
Let $V$ be a finite-dimensional irreducible  $\frak g$-module. For every element $\lambda \in \mathcal{ \frak h}^*$, we define
$$V_{\lambda} := \{ v \in V :\rho(h) = \lambda(h)v \quad \mbox{for all} \quad h \in \mathcal{\frak h}\},$$
then we have $V =   \oplus_{\lambda \in \mathcal{\frak h}^*}V_{\lambda}.$
The character of $V$   is  
 $ \ch (V) := \sum_{\lambda \in  \mathcal{\frak h}^*} (\dim V_{\lambda}) e^{\lambda}.$  
 The following formula for the character of typical irreducible modules is due to Kac \cite{Kac3}: 
\begin{equation}\label{dl0} \ch (V) =  \frac{L_1}{L_0} \sum_{w \in S_m\times S_n} \mathrm{sign} (w) e^{w(\lambda + \rho )},\end{equation}
 with $  L_1 = \sum _{\alpha \in \Delta _1^+}( e^{\alpha /2} + e^{-\alpha /2}), L_0 = \sum_{\beta \in \Delta _0^+} (e^{\beta /2} - e^{-\beta /2}).$
 
In  \cite{Zhang2},  Su and    Zhang gave an character formula 
for all finite dimension irreducible representations with any typical and atypical dominant integral weight   $\lambda $.
 The formula is too complicated to recall here, but see below for a special case.
 
\subsection{Characters of  irreducible representations of  $\frak {gl}(3|1)$}
 In this section,    we will recall formulas for the character of all typical and atypical finite-dimensional irreducible representations of $\frak {gl}(3|1)$. According to \cite[Theorem 4.9]{Zhang2}.

In $\frak {gl}(3|1)$, we have    $ \Delta _1^{+} = \{ \epsilon _1 - \epsilon _4,  \epsilon _2 - \epsilon _4, \epsilon _3 -\epsilon _4 \}$, 
   $\Delta _0^+=\{ \epsilon _1 - \epsilon _2, \epsilon _1 - \epsilon _3, \epsilon _2 - \epsilon _3\}$.  $\rho = ( \frac{1}{2}, \frac{-1}{2},
\frac{-3}{2}| \frac{-3}{2}). $ 

Set 
  $ x_1 : = e^{\epsilon _1}, x_2 := e^{\epsilon _2}, x_3 := e^{\epsilon _3}, y := e^{\epsilon _4},$ 
  $R:= (x_1 +y)(x_2+y)(x_3+y), \Pi := (x_1-x_2)(x_2-x_3)(x_1-x_3),$ 
$$a (t,u,v) := \det \left( \begin{matrix}
x_1^{t+2}&x_1^{u+1}&x_1^v\\
x_2^{t+2}&x_2^{u+1}&x_2^v\\
x_3^{t+2}&x_3^{u+1}&x_3^v
\end{matrix}\right ).$$ 
Let $ \lambda = (\lambda _1, \lambda _2,\lambda _3|\lambda _4)$
be a typical dominant integral weight. According to the character formula (\ref{dl0}), we have 
 \begin{equation*}
 \ch (V(\lambda)) = \frac{R(x_1x_2x_3)^{\lambda _3 -1}}{\Pi y^{\lambda _4}}a(\lambda _1 - \lambda _3 ,\lambda _2 -\lambda _3 ,0). 
  \end{equation*}
 
Let  $\lambda$ be an atypical weight. Then there are three possibilities:

  If $\lambda _1 + 2 = \lambda _4$, then
 \begin{multline}\label{ctdt2}
\ch (V(\lambda ))  = \frac{R}{\Pi y^{\lambda _4}}\Big[ \frac{x_1^{\lambda _1 +2}}{x_1 +y}
 ( x_2^{\lambda _2}x_3^{\lambda _3-1} - x_2^{\lambda _3 -1}
 x_3^{\lambda _2})  + \frac{x_2^{\lambda _1 +2}}{x_2 +y}  ( x_3^{\lambda _2}x_1^{\lambda _3-1} - x_3^{\lambda _3 -1}
 x_1^{\lambda _2})  \\+ 
\frac{x_3^{\lambda _1 + 2}}{x_3 + y}( x_1^{\lambda _2}x_2^{\lambda _3-1} - x_1^{\lambda _3 -1}
 x_2^{\lambda _2}) \Big].
 \end{multline}
 
If  $\lambda _2 + 1 = \lambda _4$, then
  \begin{multline}\label{ctdt3}
 \ch (V(\lambda )) = \frac{R}{\Pi y^{\lambda _4}} \Big[ \frac{x_1^{\lambda _2 +1}}{x_1 +y}
 ( x_2^{\lambda _3-1}x_3^{\lambda _1+1} - x_2^{\lambda _1+1}
 x_3^{\lambda_3-1})  + \frac{x_2^{\lambda _2 +1}}{x_2 +y}  ( x_3^{\lambda _3-1}x_1^{\lambda _1+1} - x_3^{\lambda _1+1}
 x_1^{\lambda _3-1})  \\+ 
\frac{x_3^{\lambda _2 + 1}}{x_3 + y}( x_1^{\lambda _3-1}x_2^{\lambda _1+1} - x_1^{\lambda _1+1}
 x_2^{\lambda _3-1}) \Big].
 \end{multline}
 
If   $ \lambda _3 = \lambda _4$, then
  \begin{multline}\label{ctdt4}
  \ch (V(\lambda )) = \frac{R}{\Pi y^{\lambda _4}} \Big[ \frac{x_1^{\lambda _3}}{x_1 +y}
 ( x_2^{\lambda _1+1}x_3^{\lambda _2} - x_2^{\lambda _2}
 x_3^{\lambda_1+1})  + \frac{x_2^{\lambda _3}}{x_2 +y}  ( x_3^{\lambda _1+1}x_1^{\lambda _2} - x_3^{\lambda _2}
 x_1^{\lambda _1+1})  \\+ 
\frac{x_3^{\lambda _3}}{x_3 + y}( x_1^{\lambda _1+1}x_2^{\lambda _2} - x_1^{\lambda _2}
 x_2^{\lambda _1+1}) \Big].
 \end{multline}
 \section{Double Koszul complexes}
  \subsection{The Koszul complex K}
In \cite{Manin1}  Manin suggested the following construction  to define the super determinant of  a super matrix.
Let $V^*$ denote the vector space dual to $V$ with the dual basic $\xi ^1, \xi^2,\ldots ,\xi^d$, $\xi ^i(x_j) =\delta ^i_j$. 
 The complex $K$ has its $(k,l)$-term given by $K^{k,l} := \Lambda _k\otimes S_l^*$, where $\Lambda_k$ is the $k$-th homogeneous component of
the exterior   tensor algebra over $V$,  $S_l^*$ is the $l$-th homogeneous component of the symmetric tensor algebra over $V^*$.
 The differential $d_{k,l}: K^{k,l} \longrightarrow K^{k+1,l+1}$ is given by 
\begin{equation}\label{ct1}
d_{k,l} (h\otimes \varphi ) = \sum_i h\wedge x_i \otimes \xi ^i\cdot\varphi.
\end{equation}
In fact, the construction above gives a series of complexes $K_a$:
$$\xymatrix{K_a:\cdots\ar[r]^d&\Lambda_k\otimes S_{k-a}^*\ar[r]^d&
\Lambda_{k+1}\otimes S_{k-a+1}^*\ar[r]^d&\cdots}$$
here for $k<0$ we define $\Lambda_k$ and $S_k$ to be 0. Thus each
complex $K_a$ is bounded from below.
 
It is easy to check that $d_{k,l} $ is $\frak {gl}(V)$-equivariant; hence the homology groups of this complex are representations of $\frak {gl}(V)$. On 
the other hand, one can show that the complex $(K_a,d)$  is exact everywhere if
$a\neq m-n$, and the complex $(K_{m-n},d)$ is exact everywhere except at the  term $\Lambda_m\otimes S_n^*$, 
where the homology group is one-dimensional.
This homology group defines a one-dimensional representation of $\frak {gl}(V)$. It turns out that elements of $\frak {gl}(V)$
 act on this representation by means of its super determinant.

Notice that there is another differential  
$\partial _{k,l} :  K^{k+1,l+1}\longrightarrow K^{k,l}$, which is defined as follows:
$$\xymatrix @C=2em{\partial _{k,l}: \Lambda_{k+1}\otimes S_{l+1}\ar@{^(->}[r] &V^{\otimes k+1}\otimes V^{\otimes l+1}
\ar[rrrr]^{(\Id \otimes \ev _V\otimes 
\Id)\circ (\Id \otimes  \tau _{V\otimes V^*}\otimes \Id)}&&&&V^{\otimes k} \otimes V^{*\otimes l}\ar[r]^{Y_k\otimes X_l^*}&\Lambda_k\otimes S_l^*},$$
where $$  X_n := \frac {1}{n !}\sum _{w \in \sigma _n}T_w, Y_n := \frac {1}{n !}\sum_{w \in \sigma _n}(-1)^{l(w)}T_w,$$
$$ \tau (a\otimes \varphi ) = (-1)^{\hat a.\hat \varphi }\varphi \otimes a; \text { and }   \ev(\varphi \otimes a) = \varphi (a), \mbox { where } a \in V, \varphi \in V^*;
a,\varphi - \mbox {homogeneous}.$$
  One checks that  on $ K^{k,l}$
\begin{equation}\label{ct3}
  lk d_{k-1,d-1}\partial_{k-1,l-1} +(l+1)(k+1)\partial_{k,l} d_{k,l}   =  ( l-k -n+m )\id. 
  \end{equation}
  Since $(K_\bullet , d)$ is exact, 
   $(K_\bullet ,  \partial )$   is also exact.

\subsection{The Koszul Complex $L$} There is another Koszul complex associated to $V$.
This complex was first defined by Priddy as a free resolution of $k$ as a module over the symmetric tensor
algebra of $V$ (see \cite{Manin}). As in the case of the complex $K$, the complex $L$ with $L^{p,r} := S_p\otimes \Lambda_r$ is defined as 
a series of complexes $L_a$,
$$\xymatrix{L_a:\cdots\ar[r]^{  P}&S_p\otimes \Lambda_{a-p}\ar[r]^{  P}&
S_{p-1}\otimes \Lambda_{a-p+1}\ar[r]^{  P}&\cdots}$$
with  differential $ P_{p,r}: L^{p,r} \longrightarrow L^{p-1,r+1}$  
 given by
$$ \xymatrix @C=4em{    P_{p,r}: S_p\otimes  \Lambda_r \ar@{^(->}[r] &V^{\otimes p}\otimes  V^{\otimes r}  
\ar[r]^{  X_{p-1} \otimes Y_{r+1}}&S_{p-1}\otimes  \Lambda_{r+1}}.$$
The complexes $(L_\bullet , P)$ are exact, except for $a = 0$.
\\
We also have another differential $Q_{p,r}: L^{p-1,r+1}\longrightarrow L^{p,r} $,   given by
$$ \xymatrix@C=4em { Q_{p,r}: S_{p-1}\otimes \Lambda_{r+1}\ar @{^(->}[r]&V^{\otimes p-1}\otimes V^{\otimes r+1}= V^{\otimes p}\otimes V^{\otimes r}
   \ar[r]^{\qquad \qquad X_p \otimes Y_r}& S_p\otimes \Lambda_r}.$$
 One checks that on $ L^{p,r}$ 
 \begin{equation}\label{ct60}
r (p+1)    PQ + p(r+1)Q    P   = (p+r)\Id.
\end{equation}  
Consequently the  complexes $(L_\bullet ,Q )$ are  exact too.
 
\subsection{The double Koszul complex}
The main observation of this work is the fact that the two Koszul
complexes mentioned in the previous section can be combined
into a double complex which we call the double Koszul complex.
In this section we describe this complex. An application to the
construction of irreducible representations of the super Lie algebra
$\frak{gl}(3|1)$ will be given in the next section.

For simplicity we shall  use the dot ``$ \cdot $'' to denote the tensor product.  Fix  an integer $a\geq 1$.
 We arrange the Koszul
complexes $K_{-a}, K_{-a-1}, K_{-a-2}, \ldots$  as in the diagram below. 
{\small $$ \xymatrix@R=1em{ K_{-a:}0\ar[r]&S^*_{a}\ar[r]^{d_{0,a}}&\Lambda_1\cdot S_{a+1}^*\ar[r]^{d_{1,a+1}}&\Lambda _2\cdot S_{a+2}^*
\ar[r]^{d_{2,a+2}}
&\Lambda _3\cdot S_{a+3}^*\ar[r]&\ldots\\
K_{-a-1}: &0\ar[r]&S_{a+1}^*\ar[r]^{d_{0,a+1}}&\Lambda_1\cdot S_{a+2}^* \ar[r]^{d_{1,a+2}}&\Lambda _2\cdot S_{a+3}^*\ar[r]&\ldots\\
K_{-a -2}: &&0\ar[r]&S_{a+2}^*\ar[r]^{d_{0,a+2}}&\Lambda_1\cdot S_{a+3}^* \ar[r]&\ldots }$$}
To get the entries on a column into a complex we tensor  each complex $K_{-a-i}$ with $S_{i}$,
i.e. the complex $K_{-a-1}$ is tensored with $S_1$, the complex $K_{-a-2}$ is tensored with $S_2$, etc.
Then each  column can be interpreted as the complexes
$L_j$ tensored with $S_{a+j}^*$.  
 Thus we have the following diagram where all rows are the Koszul complex $K_\bullet$ tensored with
$S_\bullet$ and the columns are the Koszul complex $L_\bullet$ tensored with $S_\bullet^*$:
\begin{equation}\label{phuckep0}{\small 
\xymatrix@R=1em{& 0&0&0&0\\
0\ar@<.5ex>[r]    &\ar[u]S_{a}^* \ar@<.5ex>[r]^d   
&\ar[u]\Lambda_1\cdot S_{a+1}^* \ar@<.5ex>[r]^d     &\ar[u]\Lambda _2\cdot S_{a+2}^* 
\ar@<.5ex>[r]^d  
  &\ar[u]\Lambda _3\cdot S_{a+3}^* \ar@<.5ex>[r]^d      &
 \ldots\\
   &0\ar@<.5ex>[r]  \ar@<1ex>[u] &S_1\cdot S_{a+1}^* \ar@<.5ex>[r]^d  
\ar@<1ex>[u]^{ P}  
&S_1\cdot \Lambda _1\cdot S_{a+2}^* \ar@<.5ex>[r]^d   \ar@<1ex>[u]^{ P} &S_1\cdot \Lambda _2\cdot S_{a+3}^*
\ar@<.5ex>[r]^d  
\ar@<1ex>[u]^{  P  }
&
  \ldots \\
 &&0\ar@<.5ex>[r]  \ar@<1ex>[u] &
  S_2\cdot S_{a+2}^* \ar@<.5ex>[r]^d    \ar@<1ex>[u]^{  P} 
&S_2\cdot \Lambda _1\cdot S_{a+3}^*\ar@<.5ex>[r]^d  
\ar@<1ex>[u]^{  P}   & 
  \ldots \\
&&&0\ar[u]&\vdots\ar[u]}}\end{equation}
  A general square in  diagram (\ref{phuckep0}) has the form
\begin{equation}\label{cell}{\small\xymatrix@R=1em{
S_i\cdot \Lambda_k\cdot S_l^*\ar[r]^{\id\otimes d}& S_i\cdot \Lambda_{k+1}\cdot S^*_{l+1}\\
S_{i+1}\cdot \Lambda_{k-1}\cdot S^*_l\ar[r]_{\id\otimes d}\ar[u]^{ P\otimes\id}&
S_{i+1}\cdot \Lambda_k\cdot S_{l+1}^*\ar[u]_{  P\otimes\id} &\mbox {with} \quad l = i+k+a.\\
}}\end{equation} 

For convenience, we   denote $ d:= \Id \otimes d,  P : =  P \otimes \Id$.  
It is easy to show that  $  P d = d  P$ for all above squares.

  We also  have an exact double Koszul complex with $d,  P$ replaced by $\partial,Q$.  
 {\small
 \begin{equation}\label{phuckep1}
\xymatrix@R=1em{& 0 \ar@{ -->}[d]  &0\ar@{ -->}[d]&0\ar@{-->}[d]&0\ar@{-->}[d]&
\\
 0  \ar@<-.5ex>@{<--}[r]_{\partial} &S_{a}^* \ar@{ -->}[d]^Q   \ar@<-.5ex>@{<--}[r]_{\partial}
& \Lambda_1\cdot S_{a+1}^* \ar@{ -->}[d]^Q\ar@<-.5ex>@{<--}[r]_{\partial}     & \Lambda _2\cdot S_{a+2}^* \ar@{ -->}[d]^Q
\ar@<-.5ex>@{<--}[r]_{\partial}
  & \Lambda _3\cdot S_{a+3}^*\ar@<-.5ex>@{<--}[r]_{\partial}  \ar@{ -->}[d]^Q    & 
 \ldots\\
   &0\ar@<-.5ex>@{<--}[r]_{\partial}  &S_1\cdot S_{a+1}^* \ar@<-.5ex>@{<--}[r]_{\partial}
\ar@{ -->}[d]^Q
&S_1\cdot \Lambda _1\cdot S_{a+2}^*\ar@<-.5ex>@{<--}[r]_{\partial}  \ar@{ -->}[d]^Q &S_1\cdot \Lambda _2\cdot S_{a+3}^*
\ar@<-.5ex>@{<--}[r]_{\partial} 
\ar@{ -->}[d]^Q
& 
  \ldots \\
 &&0\ar@<-.5ex>@{<--}[r]_{\partial} &\ar@{ -->}[d]
  S_2\cdot S_{a+2}^*\ar@<-.5ex>@{<--}[r]_{\partial}    
&\ar@{ -->}[d] S_2\cdot \Lambda _1\cdot S_{a+3}^*\ar@<-.5ex>@{<--}[r]_{\partial} 
 &  
  \ldots\\ &&&0 &\vdots }
  \end{equation}}
  The commutativity of this diagram is easy to  check. 
  
\subsection{Some remarks on the structure of the double complex}
 In this subsection we study some maps obtained from the differentials
 of the double Koszul complex.  From now, we only consider the case $(m|n) = (3|1)$.

We put the two diagrams \eqref{phuckep0} and
\eqref{phuckep1} into one: 
{\small
   \begin{equation}\label{phuckep3}
   \xymatrix{S_{i-1}\cdot S_{a+i-1}^*\ar@<.5ex>@{^(->}[r]^{d_{0,a+i-1}}   \ar@<-.5ex>@{<<--}[r]_{\partial_{0,a+i-1}}&
S_{i-1}\cdot \Lambda_1\cdot S_{a+i}^* \ar@{-->>}[d]^Q\ar@<.5ex>[r]^{d_{1,{a+i}}}   \ar@<-.5ex>@{<--} [r]_{\partial_{1,{a+i}}}
&S_{i-1}\cdot \Lambda_2 \cdot S_{{a+i}+1^*}\ar@{-->}[d]^Q\ar@<.5ex> [r]^{d_{2,{a+i}+1}}   \ar@<-.5ex>@{<--} [r]_{\partial_{2,{a+i}+1}}&\cdots\\
          &S_i\cdot S_{a+i}^* \ar@<1ex>@{^(->}[u]^{  P}\ar@<.5ex>@{^(->}[r]^{d_{0,{a+i}}}   \ar@<-.5ex>@{<<--}[r]_{\partial_{0,{a+i}}}& 
 S_i\cdot \Lambda _1.S_{{a+i}+1}^*  
\ar@<1ex>[u]^{  P}\ar@<.5ex>[r]^{d_{1,{a+i}+1}}  \ar@<-.5ex>@{<--}[r]_{\partial_{1,{a+i}+1}}   \ar@{-->>}[d]^Q&
S_i\cdot \Lambda _2\cdot S_{{a+i}+2}^*   \ar@{-->}[d]^Q\\
           & & S_{i+1}\cdot S_{{a+i}+1}^* \ar@<1ex>@{^(->}[u]^{  P} 
\ar@<.5ex>@{^(->}[r]^{d_{0,{a+i}+1}}  \ar@<-.5ex>@{<<--}[r]_{\partial_{0,{a+i}+1}}\ar@<1ex>[u]^{ P}
& S_{i+1}\cdot \Lambda _1\cdot S_{{a+i}+2}^*  \ar@<1ex>[u]^{  P}}
\end{equation}  }
       \begin{pro}\label{mde1}
The composed map   $  \partial    PQd  : S_i\cdot S_{{a+i}}^* \longrightarrow S_i\cdot S_{{a+i}}^* $ in   diagram  (\ref{phuckep3}) is an  isomorphism for all
 $i\geq 0 $. Consequently $S_i\cdot S_{a+i}^*$ is isomorphic to a direct summand of $S_{i+1}\cdot S_{{a+i}+1}^*$.
   \end{pro} 
      {\it Proof. }According to     formulas (\ref{ct3}) and (\ref{ct60}) and the commutativity between $d,P$ and $\partial, Q$, we have  
\begin{align*}
  \partial   PQd &=   \partial d - \frac{2i}{i+1}\partial Q Pd\\
 &  = \partial d - \frac{i}{i+1}QP + \frac {i(a+i)}{(i+1)({a+i}+1)}Qd\partial P\\
&=\Big [\frac{({a+i}+2)}{({a+i}+1)} - \frac { i}{i+1}\Big ]\Id + \frac {i(a+i)}{(i+1)({a+i}+1)} Qd\partial P.
\end{align*}

 We will use induction on $i$ to prove that $\partial   PQd : S_{i}\cdot S_{{a+i}}^* \longrightarrow  S_{i}\cdot S_{{a+i}}^*$ 
is diagonalizable with the set of  eigenvalues 
\begin{align*}\label{ctmd1}
 A_i :=&\Big \{\frac{(a+i+3-j)j}{(i+1)(a+i+1)}, j=1,2,\ldots,i+1
\Big \}.
\end{align*} 
For  $i=0$   the claim follows from the equation above.
 Assume that  the proposition is true for  $i-1$.

 By assumption
$  \partial PQd : S_{i-1}\cdot S_{{a+i}-1}^*\longrightarrow S_{i-1}\cdot S_{{a+i}-1}^* $ is diagonalizable with the set of 
 eigenvalues   is $A_{i-1},$
  hence $Qd\partial P: S_i\cdot S_{a+i}^*\longrightarrow  S_i\cdot S_{a+i}^*$ is diagonalizable with the set of eigenvalues is $A_{i-1}\cup \{ 0\}.$
Thus it is easy to see that   $ \partial   PQd: S_{i}\cdot S_{{a+i}}^*\longrightarrow S_{i}\cdot S_{{a+i}}^* $ 
is diagonalizable with $A_i$ the set of eigenvalues.\eee

 Consider the diagram in (\ref{phuckep0}) as an exact sequence of horizontal complexes (except for the first column) and split it into short exact sequences.   
 
{\small
\begin{equation}\label{phuckep4}
  \xymatrix{
 \ldots \ar[r] & \Ker P_{i,k}\cdot S^*_{i+k+a}\ar@<.5ex>[r]^{d'_{k,i+k+a}}\ar[d]^Q&  \Ker P_{i,k+1}\cdot S^*_{i+k+a+1}\ar[d]^Q
\ar@<.5ex>[r]^{d'_{k+1,i+k+a+1}}&  \Ker P_{i,k+2}\cdot S^*_{i+k+a+2} \ar[r]\ar[d]^Q
&\ldots \\
\ldots\ar[r]&S_{i+1}\cdot \Lambda _{k-1} \cdot S^*_{i+k+a}\ar[d]^Q\ar@<1ex>@{->>}[u]^{  P_{i+1,k-1}} 
  \ar@<.5ex>[r]^{d_{k-1,i+k+a}}
 &S_{i+1}\cdot \Lambda _k\cdot S^*_{i+k+a+1}\ar[d]^Q
\ar@<1ex>@{->>}[u]^{ P_{i+1,k}}   \ar@<.5ex>[r]^{d_{k,i+k+a+1}}
 &S_{i+1}\cdot \Lambda _{k+1}\cdot S^*_{i+k+a+2}\ar@<1ex>@{->>}[u]^{  P_{i+1,k+1}} \ar[r]\ar[d]^Q&\ldots\\
\ldots\ar[r]&\Ker P_{i+1,k-1}\cdot S^*_{i+k+a} \ar@<.5ex>[r]^{d'_{k-1,i+k+a}} \ar@<1ex> @{^(->}[u]^{i}&\Ker P_{i+1,k }
\cdot S_{i+k+a+1}^* 
\ar@<.5ex>[r]^{d'_{k,i+k+a+1}}  \ar@<1ex> @{^(->}[u]^{i}&
\Ker P_{i+1,k+1}\cdot S^*_{i+k+a+2}\ar@<1ex>@{^(->}[u]^{i}\ar[r]&\ldots }
\end{equation}}
where $d'_{k,i+k+a} := d_{k,i+k+a}\big |_{\Ker P_{i,k}\cdot S_{i+k+a}^*}$ , $\Ker P_{i,j} = \Im P_{i+1,j-1}$ for all $i\geq 0.$
 The differentials $\partial$   however do not restrict to differentials on the first and the
third horizontal complexes. 
    Consider the following part of \eqref{phuckep4}  for $i, k\geq 1$:
{\small
\begin{equation}\label{sodo}\xymatrix{ S_{i-1}\cdot  \Lambda_{k+1}\cdot S_{a+i+k}^*\ar[r] \ar @<-.5 ex>@{-->}[d]
&S_{i-1}\cdot \Lambda_{k+2}\cdot S_{a+i+k+1}\ar@<-.5ex>@{-->}[l] \ar @<-.5 ex>@{-->}[d] &\\
  \Ker P_{i-,k+1}\cdot S^*_{a+i+k}\ar[r] \ar @{^(->}[u] \ar @<-.5ex>@{-->}[d]_Q&
\Ker P_{i-1,k+2}\cdot S_{a+i+k+1}^*\ar @{^(->}[u] \ar @<-.5ex>@{-->}[d] \ar@<-.5ex>@{-->}[l] &\\
  S_i\cdot \Lambda_k \cdot S_{a+i+k}^*\ar[r]_d \ar @{->>} [u]_P &\ar@<-.5ex>@{-->}[l]_\partial \ar @{->>}[u]S_i\cdot \Lambda_{k+1} \cdot S_{a+i+k+1}^*\ar[r] \ar @<-.5 ex>@{-->}[d]& \cdot S_i \cdot \Lambda_{k+2}\cdot S^*_{a+i+k+2}
\ar@<-.5ex>@{-->}[l] \ar @<-.5 ex>@{-->}[d] \\
&  \Ker P_{i,k+1}\cdot S^*_{a+i+k+1}\ar[r] \ar @{^(->}[u] \ar @<-.5ex>@{-->}[d]_Q&
\Ker P_{i,k+2}\cdot S_{a+i+k+2}^*\ar @{^(->}[u] \ar @<-.5ex>@{-->}[d] \ar@<-.5ex>@{-->}[l]   \\
&  S_{i+1}\cdot \Lambda_{k}\cdot S^*_{a+i+k+1}\ar[r]_d  \ar  @{->>}[u]_{ P}&
S_{i+1}\cdot \Lambda_{k+1}\cdot S_{a+i+k+2}^*  \ar @{->>}[u] \ar@<-.5ex>@{-->}[l]_{\partial} }
\end{equation}  }

 \begin{pro}\label{dl1}
The composed map $$ P\partial dQ: \Ker P_{i,k+1}\cdot S_{a+i +k+1}^*\longrightarrow \Ker P_{i,k+1}\cdot S_{a+i+k+1}^*$$ 
(for $i\geq 0, k\geq 1$) in  the  diagram (\ref{sodo})
 is an  isomorphism.  Consequently $\Ker P_{i,k+1}\cdot S_{a+i+k+1}^*$ is isomorphic to a direct summand of $S_{i+1}\cdot \Im d_{k,a+i+k+1}.$
\end{pro}

{\it Proof.}
By using the method of  induction,
we will prove that $$ P\partial dQ: \Ker P_{i,k+1}\cdot S_{a+i +k+1}^*\longrightarrow \Ker P_{i,k+1}\cdot S_{a+i+k+1}^*,$$
is diagonalizable with  the set of  eigenvalues is  
{\small
\begin{align*}
&A_i := \Big \{\frac{(a+k+2i+4-j)j}{(i+1)(k+1)^2(a+i+k+2)}, j=1,2,\ldots,i+1,i+k+1
 \Big \}.
\end{align*}}

 For $i =0$, consider  the following part of \eqref{sodo}:
 {\small
\begin{equation*} \xymatrix{\Lambda_{k }\cdot   S_{a+k }^*\ar[r]_d \ar @<-.5 ex>@{-->}[d]_Q& \ar@<-.5ex>@{-->}[l]_{\partial}
  \Lambda_{k+1}\cdot   S_{a+k+1}^*\ar[r]_d \ar@<-.5ex>@{-->}[d]_Q& 
\cdot \Lambda_{k+2}\cdot S^*_{a+k+2}\ar@<-.5ex>@{-->}[d] \ar@<-.5ex>@{-->}[l]_{\partial} \\
\cdots   S_1\cdot\Lambda_{k-1}\cdot S^*_{a+k}\ar[r]_d \ar@{->>}[u]_{  P} & S_{1}\cdot \Lambda_{k}\cdot S^*_{a+k+1}\ar[r]_d \ar@{->>}[u]_{  P}\ar@<-.5ex>@{-->}[l]_{\partial}&
S_{1}\cdot \Lambda_{k+1}\cdot S_{a+k+2}^*\ar@{->>}[u]   \ar@<-.5ex>@{-->}[l]_{\partial} 
}\end{equation*} }
The composed map $ P\partial DQ: \Lambda_{k+1}\cdot S_{a+k+1}^*\longrightarrow \Lambda_{k+1}\cdot S_{a+k+1}^*.$
By means of   formulas (\ref{ct3}) and (\ref{ct60}) we have 
\begin{align*}
P\partial dQ = P\frac{[(a+3)- k(a+k+1)d\partial ]}{(k+1)(a+k+2)} Q  = \frac {(a+3)} {(k+1)(a+k+2)}\Id   - \frac {k(a+k+1)}{(k+1)(a+k+2)}d \partial.\\
\end{align*}
We have  $d\partial $  is diagonalizable with eigenvalues    $0$ and $ \frac {a+2}{(k+1)(a+k+1)}$, hence
  $ P\partial dQ$ is diagonalizable with the set of   eigenvalues    
$$ A_0 := \Big \{\frac {(a+3)}{(k+1)(a+k+2)}, 
 \frac {(a+k+3)}{(k+1)^2(a+k+2)}\Big \}.$$

For $i=1$, 
 consider in \eqref{sodo} the map $$ P\partial dQ: \Ker P_{1,k+1}\cdot S_{a+k+2}^*\longrightarrow \Ker P_{1,k+1}\cdot S_{a+k+2}^*.$$ 
On $S_i\cdot \Lambda_{k+1}\cdot S^*_{a+i+k+1}$, we have
{\small
\begin{align*}
 P\partial dQ &= P\frac{[(a+4) - k(a+k+2)d\partial ]}{(k+1)(a+k+3)}Q = \frac {(a+4)}{(k+1)(a+k+3)} PQ -\frac{k(a+k+2)}{(k+1)(a+k+3)} P d\partial Q. \end{align*}}
We have $Pd\partial Q=dPQ\partial$ and this operator can be restricted
to $\Ker P_{1,k+1}\cdot S_{a+k+2}^*$. We compute the eigenvalue of
this operator. First we have
\begin{align*}
dPQ\partial&=    d[\frac{(k+1) - (k+1) QP}{2 k}]\partial 
 =   \frac { k+1 }{2 k }d\partial  - \frac { k+1 }{2 k }dQP\partial.  
\end{align*}
Notice that $dQP\partial$ is an endomorphism of $\Im d\subset S_i\cdot \Lambda_{k+1}\cdot S^*_{a+i+k+1}$, on this space $d\partial$
operates by multiplication with $ \frac {a+3}{(k+1)(a+k+2)}$.
On the other hand from above we know the eigenvalues of
$P\partial dQ$ form the set $A_0$. Thus  
$dQP\partial$ is diagonalizable with eigenvalues $A_0\cup \{0 \}$.
Consequently $dPQ\partial$ is diagonalizable with eigenvalues
$$\Big\{\frac{a+3}{2k(a+k+2)}, \frac{a+2}{2(k+1)(a+k+2)},0\Big\}.$$
On the other hand, the restriction of $PQ$ to $\Ker P_{1,k+1}\cdot S_{a+k+2}^*$ is the multiplication with $\frac{k+2}{2(k+1)}$.
Therefore, the eigenvalues of $P\partial dQ$ are
 \begin{align*}
A_1 := &  \Big \{\frac {(a+4)(k+2)}{2(k+1)^2(a+k+3)}, \frac{(a+k+5)}{2(k+1)^2(a+k+3)}, \frac {2(a+k+4)}{2(k+1)^2(a+k+3)}
 \Big \}. 
\end{align*}

  In general, 
 we consider the composed map $$P\partial dQ: \Ker P_{i,k+1}\cdot S_{a+i+k+1}^* \longrightarrow \Ker P_{i,k+1}\cdot S_{a+i+k+1}^*
\mbox { in the  diagram (\ref{sodo}).}$$
 We have
 {\small
\begin{align*}
&P\partial dQ = P[\frac{(a+i+3) - k(a+i+k+1)d\partial }{(k+1)(a+i+k+2)}]Q\\
&= \frac {(a+i+3)}{(k+1)(a+i+k+2)}PQ- \frac{k(a+i+k+1)}{(k+1)(a+i+k+2)}dPQ\partial \\
&= \frac{(a+i+3)(i+k+1)}{(k+1)^2(i+1)(a+i+k+2)}\Id -\frac{k(a+i+k+1)}{(k+1)(a+i+k+2)}d[\frac {(i+k) - i(k+1)QP}{k(i+1)}]\partial \\
&=  \frac{(a+i+3)(i+k+1)}{(k+1)^2(i+1)(a+i+k+2)}\Id  - \frac{(i+k)(a+i+k+1)}{(k+1)(i+1)(a+i+k+2)}d\partial 
+ \frac{i(a+i+k+1)}{(i+1)(a+i+k+2)}dQP\partial .
\end{align*}}
Similar arguments shows that  $dQP\partial: \Ker P_{i,k+1}\cdot S_{a+i+k+1}^* \longrightarrow 
 \Ker P_{i,k+1}\cdot S_{a+i+k+1}^*$  is diagonalizable with the set of  eigenvalues being $A_{i-1}\cup \{0\}$. Thus  
the composed map  $  P\partial dQ: \Ker P_{i,k+1}\cdot S_{a+i+k+1}^* \longrightarrow 
 \Ker P_{i,k+1}\cdot S_{a+i+k+1}^*$  is diagonalizable with $A_i$ the set of   eigenvalues, hence is an isomorphism.\eee


    \section{Construction of irreducible representations of $\frak {gl}(V)$.}
 Let $V$ be a super vector space with super-dimension $(3|1)$. 
 In this section, using the double Koszul complex, we will construct all irreducible representations of this super algebra. To show the representations 
obtained are in fact irreducible we compute their characters.
\subsection{Combinatorial construction of irreducible representations of $\frak {gl}(V)$}
In this section, we will compute the character of the duals of irreducible direct summand of the power of the
fundamental representation $V$.
 By  the combinatorial method, we have  
 $$V^{\otimes k} = \bigoplus_{\lambda  \in \Gamma_{3,1} } I_{\lambda }^{\oplus  C_{\lambda }},$$
 where $ I_{\lambda }$ are simple, and
  $\Gamma_{3,1} $ is the set of   partitions with 
$\lambda _4 \leq 1.$  
 Since the character of $V$ is $x_1+x_2+x_3-y$, using the determinant formula (3.5) of \cite{Macdonald}, we can  compute the
character of  $I_\lambda $ for all $ \lambda \in \Gamma_{3,1}$.
 
 If  $\lambda \in \Gamma _{3,1}$ and $\lambda_3\geq 1$, we have
 \begin{equation}\label{ctdt1} 
 \ch (I_{\lambda _1,\lambda _2,\lambda _3,1^\lambda _4} ) = \frac{R(x_1x_2x_3)^{\lambda _3 -1}}{\Pi y^{\lambda _4}}
a(\lambda _1 - \lambda _3  ,\lambda _2 -\lambda _3  ,0), 
 \end{equation}
 hence 
  \begin{equation*} 
 \ch (I_{\lambda _1,\lambda_2,\lambda_3,1^\lambda_4}^*)= \frac{R (x_1x_2x_3)^{-\lambda_1}}{\Pi y^{\lambda_4+3}}a(\lambda _1-\lambda _3 ,
\lambda _1-\lambda _2 ,0). 
  \end{equation*}Thus  
  $I_{\lambda _1,\lambda_2,\lambda_3,1^\lambda_4}^*$ has highest weight   $(-\lambda _3+1,-\lambda _2+1,-\lambda _1 +1|\lambda_4+3).$
  Therefore we have $$\ch (I_{\lambda _1,\lambda _2,\lambda _3,1^\lambda _4} ) = \ch (V(\lambda _1,\lambda _2,\lambda _3|-\lambda _4)), $$
  $$\ch (I_{\lambda _1,\lambda_2,\lambda_3,1^\lambda_4}^*) = \ch (V(-\lambda _3+1,-\lambda _2+1,-\lambda _1 +1|\lambda_4+3)).$$
  
  Now,
 \begin{align*} 
\ch (I_{\lambda_1,\lambda_2,0,0}) &= \frac{R}{\Pi }\Big [  \frac{x_2^{\lambda_1+1}x_3^{\lambda_2} - x_2^{\lambda_2} x_3^{\lambda_1+1}}
{x_1 + y} 
+ \frac{x_3^{\lambda_1+1}x_1^{\lambda_2} - x_3^{\lambda_2} x_1^{\lambda_1+1}}{x_2 + y}
+  \frac{x_1^{\lambda_1+1}x_2^{\lambda_2} - x_1^{\lambda_2} x_2^{\lambda_1+1}}{x_3 + y} \Big], 
\end{align*}
hence 
  \begin{multline*} 
\ch (I_{\lambda_1,\lambda_2,0,0}^*) = \frac{R}{\Pi y^2} \Big[ \frac{x_1^2}{x_1 + y}(x_2^{-\lambda_2+1}x_3^{-\lambda_1} - x_2^{-\lambda_1}
x_3^{-\lambda_2+1}) +
  \frac{x_2^2}{x_2 + y}(x_3^{-\lambda_2+1}x_1^{-\lambda_1} - x_3^{-\lambda_1}x_1^{-\lambda_2+1}) \\
+ \frac{x_3^2}{x_1 + y}(x_1^{-\lambda_2+1}x_2^{-\lambda_1} - x_1^{-\lambda_1}x_2^{-\lambda_2+1}) \Big]. 
\end{multline*} 
Therefore   $I_{\lambda_1,\lambda_2,0,0}^*$ has highest weight  
$(0,-\lambda _2+1,-\lambda _1+1|2).$
Thus we have 
$$ \ch (I_{\lambda_1,\lambda_2,0,0}) = \ch (V(\lambda_1,\lambda_2,0|0)),$$
$$   \ch (I_{\lambda_1,\lambda_2,0,0}^*) = \ch (V(0,-\lambda _2+1,-\lambda _1+1|2)).$$
Further we have
 \begin{equation}\label{ctdt8} 
\ch (I_{\lambda _1,0,0,0}) = \frac{1}{\Pi } \Big[  x_2^{\lambda _1+1}(x_2+y)(x_3-x_1) +   x_3^{\lambda _1+1}(x_3+y)(x_1-x_2) +
 x_1^{\lambda _1+1}(x_1+y)(x_2-x_3) \Big],  
 \end{equation}and  hence 
 \begin{align*} 
\ch (I_{\lambda _1,0,0,0}^*)  
&= \frac{1}{\Pi y}\Big [  x_1^2(-x_2^{-\lambda _1+1}x_3 + x_2  x_3^{-\lambda _1+1}) + x_2^2(-x_3^{-\lambda _1+1}x_1 
+ x_3  x_1^{-\lambda _1+1})\\
 &+
 x_3^2(-x_1^{-\lambda _1+1}x_2 
+ x_1  x_2^{-\lambda _1+1}) 
+ x_1^2y(-x_2^{-\lambda _1}x_3 +x_2  x_3^{-\lambda _1}) \\
&+ x^2y(-x_3^{-\lambda _1}x_1 + x_3  x_1^{-\lambda _1})
 +
 x_3^2y(-x_1^{-\lambda _1}x_2 + x_1  x_2^{-\lambda _1}) \Big ]. 
  \end{align*} 
Therefore $I_{\lambda_1,0,0,0}^*$ has highest weight  $(0,0,-\lambda _1 +1|1)$.
Thus we have
$$\ch (I_{\lambda_1,0,0,0}) = \ch (V(\lambda_1,0,0|0)), $$
$$\ch (I_{\lambda_1,0,0,0}^*) = \ch (V(0,0,-\lambda _1 +1|1)).$$
   \subsection{Construct  representations by using Koszul complex $K$}
Consider complexes  $K_a$, with  $a: = k-l  \neq 2$.
$$K_a: ....\longrightarrow \Lambda _k.S_l ^*\longrightarrow  \Lambda _{k+1}.S_{l+1}^*\longrightarrow \Lambda _{k+2}.S_{l+2}^*\longrightarrow ....,$$
 By using the exactness property of the Koszul complex $K$, we will construct a class of irreducible representations of $\frak {gl}(3|1)$. According to (\ref{ct3}) we have  
\begin{equation}\label{xdanh}
 \Lambda _k.S_l^* \cong  \Im d_{k-1,l-1} \oplus \Im d_{k,l}.
 \end{equation} 
 Consequently, we have (\cite{dh2}).
\begin{pro}
The module $\Im d_{k+1,l+1}$ is simple   for all pairs $(k,l)$ with $l,k \geq 1, k-l \neq 2.$
\end{pro}
  
  Using induction,  we find that
 \begin{equation}\label{ctdt7} 
 \ch (\Im d_{k,l}) = \frac{Ry^{k-3}}{\Pi (x_1x_2x_3)^l}a (l ,l ,0).
\end{equation}
Set $M^{m,p} := \Im d_{m+2,m+p} . H_{3,1}^{\otimes m-1}$ with  $H_{3,1}:= \Ker d_{3,1}/\Im d_{2,0}$. We have 
\begin{equation*} 
\ch (M^{m,p} )= \frac{R}{\Pi (x_1x_2x_3)^{p+1}}a (m+p , m+p , 0) = \ch (V(m,m,-p|0)).
\end{equation*} 
  Hence
  \begin{equation*} 
  \ch (M^{m,p })^* =\frac {R (x_1x_2x_3)^{-m}}{\Pi y^3}a(m+p , 0,0)= \ch V(p+1,-m+1,-m+1|3).
  \end{equation*}
   Therefore,  $M$ is isomorphic to $V(m,m,-p|0) $, and $M^*$ is isomorphic to $V(p+1,-m+1,-m+1|3). $
  
Thus every irreducible representation with highest weight in the set
 $$\{ (m,n,p|-q), (-p,-n,-m|q): (m,n,p,q) \in \Gamma _{3,1}\}\cup   
\{(m,m,-p|0) ,  (p,-m,-m|0):   m,p\geq 1\} $$
is  constructed.

It remains to construct representations with   highest weights in the set 
 $$\{(n,0,-p|0): n,p\geq 1\} \cup \{ (m+a,m, -p|,0): m,a,p\geq 1\}.$$
\subsection{Construct representations by using double Koszul complex.}
    According to Prop.  \ref {mde1}, there exists $Y$ such that 
$ S_{n}\cdot S_{p}^* = S_{n-1}\cdot S_{p-1}^* \oplus Y$. It is easy to compute for $n,p\geq 1$ that
\begin{equation}\label{ctdt11}
  \ch (Y )= \frac{(x_1x_2x_3)R}{\Pi y}\Big[  \frac{x_2^{-p-1}x_3^n - x_2^nx_3^{-p-1}}{x_1+y} +  
\frac{x_3^{-p-1}x_1^n - x_3^nx_1^{-p-1}}{x_2+y} + \frac{x_1^{-p-1}x_2^n - x_1^nx_2^{-p-1}}{x_3+y}  \Big ]. 
\end{equation}
Hence, $Y$ has highest weight   
$ (n,0,-p+1|1).$
 
Next, we will construct representations  having   highest weights   in the set  $\{ \lambda = (m+t,m,-p,0): m,p,t \geq 1\}.$
   According to Proposition \ref{dl1}, we  have $$ S_1\cdot \Im d_{2,m+1} = \Lambda _3\cdot S_{m+1}^* \oplus Z_1,$$ hence 
 \begin{equation}\label{ctdt22}
\ch (Z_1 )= \frac{R}{\Pi y(x_1x_2x_3)^{m+1}}a(m+2,m+1,0)= \ch (V(2,1,-m+1|1)),
\end{equation} 
 therefore  $Z_1$   is isomorphic to $V(2,1,-m+1|1).$
 
In general, according to Proposition \ref {dl1}, we have 
   $\Im (\Id_{I_{k,0,0,0}}\cdot d_{l,m}) = \Ker P_{k,l}\cdot S_m^* \oplus  Z_k,$  where $\Ker P_{k,l} \cong I_{k,1^l}$. Therefore 
$\ch (Z_k) =  \ch [\Im (\Id_{I_{k,0,0,0}}\cdot d_{l,m})] - \ch (I_{k,1^l}\cdot S_m^*).$

According to (\ref{ctdt1}),   (\ref{ctdt8})  and (\ref{ctdt7}), we have 
\begin{equation}\label{ctdt24}
 \ch Z_k=   \frac{R(x_1x_2x_3)^{-m}y^{l-3}}{\Pi }a(k+m,m-1,0).
  \end{equation}  Set $M := Z_k\cdot I_{1,1,1,-1}^{\otimes (l-2)}$, then 
  \begin{equation}\label{ctdt25}
 \ch (M )=     \frac{R(x_1x_2x_3)^{-p}}{\Pi y}a(m+p+t-1, m+p-1, 0) = \ch (V(m+t,m,-p+1|1)),
\end{equation} where $t :=l-1, 
p := -m-2+l.$
According to (\ref{ctdt25}), we have 
\begin{equation}\label{ctdt26}
 \ch (M^*)= \frac{R(x_1x_2x_3)^{-m-a}}{\Pi y^2}a(m+t+p-1, t,0)=\ch (V(p+1,-m+1,-m-t+1|3)).
   \end{equation}
Therefore $M$ is isomorphic to $ V(m+t,m,-p+1|1)$, $M^*$ is isomorphic to $V(p+1,-m+1,-m-t+1|3). $

Thus, for any integrable dominant  weight $\lambda = (\lambda_1, \lambda_2, \lambda_3| \lambda_4 ) $, we have constructed a representation
which has highest weight    $ \lambda $ and has character   equal to the character of  the   irreducible representation with highest weight $\lambda$.

\end{document}